\documentclass[reqno]{amsart}
\usepackage{amssymb}
\usepackage{amsfonts}
\usepackage{latexsym}

\newtheorem{theorem}{Theorem}[section]
\newtheorem{theorem/definition}{Theorem/Definition}[section]

\newtheorem{proposition}{Proposition}[section]
\newtheorem{lemma}{Lemma}[section]

\newtheorem{corollary}{Corollary}[section]
\newtheorem{conjecture}{Conjecture}[section]
\theoremstyle{remark}
\newtheorem{remark}{Remark}[section]
\theoremstyle{definition}
\newtheorem{example}{Example}[section]
\newtheorem{definition}{Definition}[section]
\newcommand{\bC}{{\mathbb C}}

\newcommand{\bR}{{\mathbb R}}
\newcommand{\bZ}{{\mathbb Z}}
\newcommand{\cL}{{\mathcal L}}
\newcommand{\cS}{{\mathcal S}}
\newcommand{\cO}{{\mathcal O}}
\newcommand{\cD}{{\mathcal D}}

\newcommand{\cF}{{\mathcal F}}
\newcommand{\cE}{{\mathcal E}}

\newcommand{\reg}{\text{reg}}
\newcommand{\GL}{\text{GL}}

\newcommand{\Pic}{\text{Pic}}
\newcommand{\Hom}{\text{Hom}}
\newcommand{\td}{\text{d}}
\newcommand{\Res}{\text{Res}}
\begin{document}

\title{On elliptic Dunkl operators}

\author{Pavel Etingof}
\address{Department of Mathematics, Massachusetts Institute of Technology,
Cambridge, MA 02139, USA}
\email{etingof@math.mit.edu}

\author{Xiaoguang Ma}
\address{Department of Mathematics, Massachusetts Institute of Technology,
Cambridge, MA 02139, USA}
\email{xma@math.mit.edu}

\maketitle

\section{Introduction}

Elliptic Dunkl operators for Weyl groups were introduced in
\cite{BFV}. Another version of such operators was considered by
Cherednik \cite{Ch1}, who used them to prove the quantum integrability of
the elliptic Calogero-Moser systems. The goal of the present 
paper is to define elliptic Dunkl operators for any finite group $W$
acting on a (compact) complex torus\footnote{We note that
although we work with a general finite group $W$, the theory
essentially reduces to the case when $W$ is a crystallographic
reflection group (\cite{GM},5.1), because $W$ can be replaced by its subgroup
generated by reflections.} $X$. We attach such a set of operators 
to any topologically trivial holomorphic line bundle ${\mathcal
L}$ on $X$ with trivial stabilizer in $W$, and any flat
holomorphic connection $\nabla$ on this bundle. 
In the case when $W$ is the Weyl group of a root system, and $X$ is the space of
homomorphisms from the root lattice to the elliptic curve,  
our operators coincide with those of \cite{BFV}. We prove that the
elliptic Dunkl operators commute, and show that the monodromy 
of the holonomic system of differential equations defined by them
gives rise to a family of $|W|$-dimensional representations of the Hecke
algebra ${\mathcal H}_\tau(X,W)$ of the orbifold $X/W$, defined in \cite{E};
conjecturally, this gives generic irreducible representations of
this algebra. In the case of Weyl groups, the algebra
${\mathcal H}_\tau(X,W)$ is the double affine
Hecke algebra (DAHA) of Cherednik \cite{Ch2}, while in the case 
$W=S_n\ltimes (\Bbb Z/\ell Z)^n$, $\ell=2,3,4,6$, 
it is the generalized DAHA introduced in \cite{EGO},
and we reproduce known families of representations of these
algebras. We also explain how to use 
the elliptic Dunkl operators to construct representations from
category ${\mathcal O}$ over the elliptic Cherednik algebra, i.e.
the Cherednik algebra of the
orbifold $X/G$ defined in \cite{E}. 

{\bf Acknowledgements.} 
The work of P.E. was  partially supported by the NSF grant
DMS-0504847.

\section{Preliminaries} 

\subsection{Finite group actions on complex tori}

Let $V$ be a finite dimensional complex vector space. 
A nontrivial element $g\in \GL(V)$ is called a reflection 
if it is semisimple and fixes a hyperplane in $V$ pointwise. 
 
Let $W$ be a finite subgroup of $GL(V)$, and $\Gamma \subset V$ is a lattice which is preserved by $W$. Then we get a $W$-action on the complex torus $X=V/\Gamma$. 
For any reflection $g\in W$, let $X^{g}$ be the set of $x\in X$
s.t. $gx=x$. A reflection hypertorus is any connected component
of $X^{g}$ which has codimension $1$. Let $X_{\reg}$ be the
complement of reflection hypertori in $X$.

Let $H$ be a reflection hypertorus.
Let $W_{H}\subset W$ be the stabilizer of a generic point in
$H$. Then  $W_{H}$ is a cyclic group with order $n_{H}$. The
generator $g_{H}$ is the element in $W_{H}$ with determinant
$\exp(2\pi{\rm i}/n_{H})$. Let $\cS$ denote the set of pairs $(H,j)$, where $H$ is a reflection hypertorus and $j=1,\ldots, n_{H}-1$.

Under the $g_{H}$-action, we have a decomposition:
\begin{equation*}\label{decomp}
V=V^{g_{H}}\oplus  {V}_{H},
\end{equation*}
where $V^{g_{H}}$ is the codimension 1 subspace of $V$ with
a trivial action of $g_H$, and $
{V}_{H}=((V^{*})^{g_H})^{\perp}$, 
which is a $g_H$-invariant $1$-dimensional space. We also have a similar
decomposition on the dual space: $V^{*}=(V^{*})^{g_{H}}\oplus  {V}_{H}^{*}$.

\subsection{Holomorphic line bundles on complex tori}
Let us recall the theory of holomorphic line bundles on complex tori (see \cite{Mu} or \cite{La} for more details). 

Let $X=V/\Gamma$ be a complex torus.  Any holomorphic line bundle $\cL$ on $X$ is a quotient of $V\times \bC$ by the $\Gamma$ action: 
 \begin{equation*}
 \gamma:(z,\xi)\mapsto(z+\gamma,\chi(\gamma,z)\xi),\end{equation*} 
 where $\chi(\gamma,\cdot):V\to \bC^{*}$ is a holomorphic function s.t. 
\begin{equation*}
\chi(\gamma_{1}+\gamma_{2},z)=\chi(\gamma_{1},z+\gamma_{2})\chi(\gamma_{2},z).
\end{equation*} 
 Denote by $\cL(\chi)$ the line bundle corresponding to $\chi$.
 
 Let $V^{\vee}=\Hom_{\bar{\bC}}(V, \bC)$ be the vector space of $\bC$-antilinear forms on $V$. We have a nondegenerate $\bR$-bilinear form 
 \begin{equation*}
 \omega:V^{\vee}\times V\to \bR, \omega(\alpha,v)=\text{Im } \alpha(v).
 \end{equation*}
 
Then we define $\Gamma^{\vee}=\{\alpha\in V^{\vee}|\omega(\alpha,\Gamma)\subset \bZ\}$. It is easy to see that $\Gamma^{\vee}$ is a lattice in $V^{\vee}$ and we have the dual torus $X^{\vee}=V^{\vee}/\Gamma^{\vee}$. 
 
To any element $\alpha\in X^{\vee}$,  we can associate a line bundle $\cL_{\alpha}=\cL(\chi_{\alpha})$, where $\chi_{\alpha}(\gamma,z)=\exp(2\pi {\rm i}\omega(\alpha,\gamma))$. This is a  topologically trivial line bundle on $X$.

 \begin{proposition}\label{pic-tor}
 The map $\alpha\mapsto \cL_{\alpha}$ is an isomorphism of groups
$$X^{\vee}\to \Pic^{0}(X).$$
 \end{proposition}

Now suppose that a finite group $W$ acts faithfully on $V$ and preserves a lattice $\Gamma$. By using the bilinear form $\omega$, we can define the dual $W$-action on $V^{\vee}$ which preserves the dual lattice $\Gamma^{\vee}$. So we have an action of $W$ on the complex torus $X$ and its dual $X^{\vee}$. 

We define a $W$-action on $\Pic^{0}(X)$ by:
 \begin{equation*}
 w: \cL_{\alpha}\mapsto \cL^{w}_{\alpha}=\cL_{w\alpha},
  \end{equation*}  
We have $(\cL^g)^h=\cL^{hg}$. 
   
\subsection{The Poincar\'e residue}

Suppose $\alpha$ is a meromorphic $1$-form on an $n$-dimensional complex
 manifold $X$ with a simple pole on a smooth hypersurface $Z\subset X$ and no other singularities. Near any point of $Z$, we can choose local coordinates $(z_{1},\ldots,z_{n})$ on $X$ s.t. $Z$ is locally defined by the equation $z_{1}=0$.
Then $\alpha$ can be locally expressed as 
\begin{equation*}
\alpha=\frac{1}{z_{1}}\sum_{i=1}^{n}\beta_{i}(z_{1},\ldots,z_{n})\td z_{i},
\end{equation*} 
where $\beta_{i}$'s are holomorphic. Then $\beta_{1}|_{Z}$ is a holomorphic function on $Z$, and it does not depend on the choice of the coordinates. 

We define the Poincar\'e residue of $\alpha$ at $Z$ to be $\Res_{Z}(\alpha)=\beta_{1}|_{Z}$.

More generally, let $\cE$ be a holomorphic vector bundle on $X$, and $s$ be a meromorphic section of $\cE\otimes T^{*}X$ which has a simple pole on a smooth hypersurface $Z\subset X$ and no other singularities. Similarly to the above, we define the Poincar\'e residue of $s$ to be an element in $\Gamma(Z,\cE|_{Z})$ denoted by $\Res_{Z}(s)$.

\section{Construction of elliptic dunkl operators}

\subsection{The sections $f_{H,j}^\cL$}
  
The goal of this subsection is to define certain 
meromorphic sections $f_{H,j}^\cL$ of the 
bundle $(\cL^{g_H^j})^*\otimes \cL\otimes V_H^*$
which are used in the definition of elliptic Dunkl operators. 

The line bundle 
$(\cL^{w}_{\alpha})^{*}\otimes\cL_{\alpha}$ is topologically
trivial, and it is holomorphically trivial if and only if $\alpha$ is a fixed point of $w$. 
Since $W$ acts faithfully on $V$, we can always find a point $\alpha\in X^{\vee}$ which is not fixed by any $w\in W$, i.e., there exists a  topologically trivial line bundle $\cL:=\cL_{\alpha}$ such that $(\cL^{w})^{*}\otimes\cL$ is nontrivial for any $w\in W$. From now on, we fix such a line bundle.

Let $H\subset X$ be a reflection hypertorus. We have the following lemma.

\begin{lemma}\label{mero-section}
For $j=1,\ldots,n_{H}-1$, the holomorphic line bundle $(\cL^{g_{H}^{j}})^{*}\otimes\cL$ has a global meromorphic section $s$ which has a simple pole on $H$ and no other singularities. Such $s$ is unique up to a scalar.
\end{lemma}

\begin{proof}
Let $\bar{H}=\{x\in X|x+H=H\}$. Then $\bar{H}$ is a complex torus.

It is sufficient to assume in the proof that  $H=\bar{H}$.

We have a short exact sequence of complex tori:
\begin{equation*}\label{fiberation}
0\to H\xrightarrow{\mu}  X\xrightarrow{\nu} E\to 0, \text{ where } E=X/H.
\end{equation*}
It induces a short exact sequence for the dual tori:
\begin{equation*}
0\to E^{\vee}\to X^{\vee}\to H^{\vee}\to 0,
\end{equation*}

which can be written using the isomorphism of Prop \ref{pic-tor}
\begin{equation*}
1\to \Pic^{0}(E)\xrightarrow{\nu^{*}} \Pic^{0}(X)\xrightarrow{\mu^{*}}  \Pic^{0}(H)\to 1.
\end{equation*}

Since $\mu^{*}((\cL^{g_{H}^{j}})^{*}\otimes\cL)$ is trivial, there exists a topologically trivial line bundle $\cL'$ on $E$ such that $\nu^{*}\cL'=(\cL^{g_{H}^{j}})^{*}\otimes\cL$.  
It is well known that $\cL'$ has a unique meromorphic section, up to
a scalar, which has simple pole at $0$. 
Then $s=\nu^{*}s'$ is the required section of the bundle $(\cL^{g_{H}^{j}})^{*}\otimes\cL$ on $X$.

Now we prove the uniqueness of $s$ up to a scalar.    

The section $s$ can be viewed as a global holomorphic section of the line bundle $\cF=(\cL^{g_{H}^{j}})^{*}\otimes\cL\otimes \cO(H)$.
Since $\cO(H)$ is the pullback of $\cO(0)$ on $E$, $\cF$ is the pullback of the bundle $\cL^{'}\otimes \cO(0)$ on $E$.  
So $H^{0}(X,\cF)\simeq H^{0}(E,\cL^{'}\otimes \cO(0))=\bC$ and $s$ is unique up to a scalar.
\end{proof} 

Now choose a nonzero element $\alpha\in  {V}_{H}^{*}$ and consider $s\otimes \alpha$, where $s$ is the global meromorphic section in Lemma \ref{mero-section}. 
Then $s\otimes \alpha$ is a global section of the bundle $(\cL^{g^{j}_{H}})^{*}\otimes \cL \otimes  {V}_{H}^{*}$. Its only singularity is a simple pole at $H$, and it is defined by this condition uniquely up to scaling.

Next, observe that since $X$ is a torus, the bundle $T^{*}X$ 
is canonically trivial, and we can canonically identify 
the fibers of $T^{*}X$ with $V^{*}$. Thus we may consider the 
Poincar\'e residue $\Res_{H}(s\otimes \alpha)$ which is an element in $\Gamma(H, ((\cL^{g^{j}_{H}})^{*}\otimes \cL)|_{H})$. 
Since $((\cL^{g^{j}_{H}})^{*}\otimes \cL)|_{H}$ is trivial, $\Res_{H}(s\otimes \alpha)$ is a holomorphic function on $H$. Since $H$ is compact, $\Res_{H}(s\otimes \alpha)$ is a constant. Then by fixing  this constant, we can fix $s\otimes \alpha$ uniquely, i.e., we have the following lemma:

\begin{lemma}\label{section}
For any $(H,j)\in \cS$, we have  a unique global meromorphic
section $f_{H,j}^\cL$ 
of the bundle $(\cL^{g^{j}_{H}})^{*}\otimes\cL\otimes
{V}_{H}^{*}$, such that it has a simple pole on $H$, no other
singularities, and  has residue 
$1$ on $H$.
\end{lemma}

\subsection{Elliptic Dunkl operators}

For any $g\in W$, we have a $W$-action on $\cS$: $g(H,j)=(gH,j)$.
Let $C$ be a $W$-invariant function on $\cS$. Choose a holomorphic flat connection $\nabla$ on $\cL$. We have the following definition:  

\begin{definition}[Elliptic Dunkl operators] 
For any $v\in V$, we define \emph{ the elliptic Dunkl operator}
corresponding to $v$ to be the following operator acting on the local meromorphic sections of $\cL$:
\begin{equation*}\label{Ell-Dunkl}
\cD^{\cL,\nabla}_{v,C}=\nabla_v-\sum_{(H,j)\in \cS}C(H,j)\langle f_{H,j}^{\cL},v\rangle g_{H}^{j},
\end{equation*}
where $\nabla_{v}$ is the covariant derivative along $v$
corresponding to the connection $\nabla$, and $\langle ,\rangle $ is
the natural pairing between $V$ and $V^{*}$.
\end{definition}

\begin{remark}
Let $\nabla,\nabla'$ be two flat holomorphic connections on
$\cL$. Then $\nabla-\nabla'=\xi$, where $\xi\in V^*$ is a
holomorphic 1-form on $X$. Therefore, 
\begin{equation*}
\cD^{\cL,\nabla}_{v,C}-\cD^{\cL,\nabla'}_{v,C}=\xi(v).
\end{equation*}
Thus elliptic Dunkl operators attached to different flat
connections on the same line bundle $\cL$ differ by additive
constants. 
\end{remark} 

For simplicity, we will use the same notation $\nabla$ 
for the connection on each bundle $\cL^w$ obtained 
from the connection $\nabla$ on
$\cL$ by the action of $w\in W$. Then we have the 
following result on the equivariance of the elliptic Dunkl
operators under the action of $W$. 

\begin{proposition}\label{Winv}
One has 
\begin{equation*}
w\circ \cD^{\cL,\nabla}_{v,C}\circ w^{-1}=
\cD^{\cL^w,\nabla}_{wv,C}.
\end{equation*}
\end{proposition}

\section{The commutativity theorem}

\subsection{The elliptic Cherednik algebra}

Let $c(H,j)=\frac{1}{2}(e^{-2\pi {\rm i}j/n_H}-1)C(H,j)$,
and set $c(H,0)=0$. Recall from \cite{E} that the 
sheaf of algebras $H_{1,c,0,X,W}$ on $X/W$ is 
defined as follows. Let $\bar U$ be a small open set in $X/W$, and 
$U$ be its preimage in $X$. Then the algebra $H_{1,c,0,X,W}(U)=H_{1,c,0}(U,W)$
is generated by the algebra of holomorphic functions
${\mathcal O}(U)$, the group $W$, and Dunkl-Opdam operators 
\begin{equation*}
D_{v,\phi}=\partial_v-\sum_{(H,j)\in {\mathcal
S}}C(H,j)\langle\phi_H,v\rangle g_H^j,
\end{equation*}
where $\phi=(\phi_H)$ is a collection of 1-forms on $U$ which 
locally near $H$ have the form $\phi_H=\td\log \ell_H+\phi_H'$,
 $\ell_H$ being a nonzero holomorphic function with a simple zero along
$H$, and $\phi_H'$ is holomorphic. 
For brevity we will denote 
this sheaf by $H_{c,X,W}$. It is called the Cherednik algebra of
the orbifold $X/W$ attached to the parameter $c$, or the {\it elliptic
Cherednik algebra.} 

The sheaf $H_{c,X,W}$ sits inside $W\ltimes D_{X_{\reg}}$, 
where $D_{X_{\reg}}$ is the sheaf of differential operators on $X$
with poles on the reflection hypertori.
Thus the sheaf $H_{c,X,W}$ has a filtration by order of differential operators.
It is known \cite{E} that we have $F^0H_{c,X,W}=W\ltimes
{\mathcal O}_X$. 
  
\subsection{The commutativity theorem}

One of the main theorems of this paper is the following result. 

\begin{theorem}\label{commu}
The elliptic Dunkl operators commute, i.e.
$[\cD^{\cL,\nabla}_{v,C},\cD^{\cL,\nabla}_{u,C}]=0$.
\end{theorem}

\begin{proof}
Since $\langle f_{H,j}^{\cL},v\rangle$ depends
only on the projection of $v$ to $ V_H$,
which is a 1-dimensional space,  
it is easy to check that the commutator
$[\cD^{\cL,\nabla}_{v,C},\cD^{\cL,\nabla}_{u,C}]$
does not have differential terms. In other words, we have 
\begin{equation*}
[\cD^{\cL,\nabla}_{v,C},\cD^{\cL,\nabla}_{u,C}]=
\sum_{g\in W}\varphi_gg
\end{equation*}
where $\varphi_g$ is a meromorphic section of the line bundle
$(\cL^g)^*\otimes \cL$.

We claim that $\varphi_1=0$. Indeed,
write $\cD^{\cL,\nabla}_{v,C}$ in the form 
\begin{equation*}
\cD^{\cL,\nabla}_{v,C}=\nabla_v-\sum_{H}\langle F_H,v\rangle,
\end{equation*}
where $F_H=\sum_{j=1}^{n_H-1}C(H,j)f_{H,j}^{\cL}g_H^j$. 
To show that $\varphi_1=0$, it suffices to show that 
\begin{equation*}
[\langle F_H,v\rangle,\langle F_K,u\rangle]+
[\langle F_K,v\rangle,\langle F_H,u\rangle]=0
\end{equation*}
if $W_H\cap W_K\ne 1$. 
But this is obvious, given that $\langle F_H,v\rangle$ depends
only on the projection of $v$ to $ V_H$, which is
1-dimensional, and $ V_H= V_K$ once $W_H\cap W_K\ne
1$.   

The rest of the proof of the theorem is based on the following key lemma. 

\begin{lemma}\label{holo}
The sections $\varphi_g$ are holomorphic. 
\end{lemma}

The lemma clearly implies the theorem, since the bundle 
$(\cL^g)^*\otimes\cL$ is a topologically, but not holomorphically, trivial 
bundle, and hence every holomorphic section of this bundle is
zero. 

\begin{proof} (of Lemma \ref{holo})
The lemma is proved by local analysis, i.e., essentially, by
reduction to the case of usual (rational) Dunkl-Opdam operators,
\cite{DO}. Namely, 
it is sufficient to show that $\varphi_g$ are regular when
restricted to a small $W$-invariant neighborhood 
$X_b$ of $Wb$, where $b\in X$ is
an arbitrary point. Let $W_b$ be the stabilizer of $b$ in $W$. 
Then $X_b$ is a union of $|W/W_b|$ small balls around the points
of the orbit $Wb$. Let us pick a trivialization of $\cL$ on
$X_b$. This trivialization defines a trivialization of the line
bundle $\cL^w$ for every $w\in W$. With these trivializations,
the elliptic Dunkl operators $\cD^{\cL,\nabla}_{v,C}$ become
operators acting on meromorphic functions on $X_b$. 
 
The remainder of the proof is based on the theory of Cherednik
algebras for orbifolds. Namely, 
it is clear from the definition of the elliptic Dunkl 
operators that they belong to the algebra $H_{c,X,W}(X_b)$. 
Since $F^0H_{c,X,W}=W\ltimes {\mathcal O}_X$, 
this implies that the sections $\varphi_g$, upon trivialization, 
become holomorphic functions on $X_b$. This proves the lemma. 
\end{proof} 

\end{proof} 

\section{Representations of elliptic Cherednik algebras 
arising from elliptic Dunkl operators}

In this section we will use elliptic Dunkl operators to construct
representations of the sheaf of elliptic Cherednik algebras
$H_{c,X,W}$ on the sheaf ${\mathcal F}:=\oplus_{w\in W}{\cL^*}^w$. 

Let us write the elliptic Dunkl
operator in the form 
\begin{equation*}
\cD^{\cL,\nabla}_{v,C}=\nabla_v-\sum_{g\in W}\langle
F_{C,g}^\cL,v\rangle g,
\end{equation*}
where 
\begin{equation*}
F_{C,g}^\cL=\sum_{(H,j): g_H^j=g}C(H,j) f_{H,j}^\cL
\end{equation*}
is a section of $(\cL^g)^*\otimes \cL\otimes
V^*$. Note that $F_{C,g}^\cL=0$ unless $g$ is a reflection. 

\begin{lemma}\label{inv}
We have:
\item{(i)\ }  ${\rm Ad}w(F_{C,g}^\cL)=F_{C,wgw^{-1}}^{\cL^w}$,
where ${\rm Ad}w$ stands for the adjoint action of $w$; 

\item{(ii)\ } $\nabla_u\langle F_{C,g}^\cL,v\rangle= \nabla_v\langle F_{C,g}^\cL,u\rangle; 
$
\item{(iii)\ }  $\sum_{h,g:hg=k}\langle F_{C,g}^\cL,v\rangle\langle
F_{C,h}^{\cL^g},gu\rangle=
\sum_{h,g:hg=k}\langle F_{C,g}^\cL,u\rangle\langle
F_{C,h}^{\cL^g},gv\rangle;$
\item{(iv)\ } $\sum_{h,g:hg=k}\langle F_{C,g}^\cL,v\rangle\langle F_{C,h}^{\cL^g},u\rangle=\sum_{h,g:hg=k}\langle F_{C,g}^\cL,u\rangle\langle F_{C,h}^{\cL^g},v\rangle.$
\end{lemma} 

\begin{proof}
Statement (i) follows from Proposition \ref{Winv}, 
Statements (ii),(iii) follow from the commutativity
of the elliptic Dunkl operators, using (i). 
Statement (iv), using (iii), reduces to the identity 
\begin{equation*}
\sum_{h,g:hg=k}\biggl(\langle F_{C,g}^\cL,v\rangle\langle
F_{C,h}^{\cL^g},u-gu\rangle-
\langle F_{C,g}^\cL,u\rangle\langle
F_{C,h}^{\cL^g},v-gv\rangle\biggr)=0.
\end{equation*}
Every summand in this sum is a skew symmetric bilinear form in $u,v$ which
factors through ${\rm Im}(1-g)$. But if $F_{C,g}^\cL$ is nonzero, then
$g$ is a reflection, and hence ${\rm Im}(1-g)$ is a 1-dimensional space. 
This means that every summand in this sum is zero, and 
the identity follows.  
\end{proof}

Now we will define the representation of the elliptic Cherednik
algebra. We start by defining an action $\rho=\rho_{\cL,\nabla}$ of the sheaf $W\ltimes
D_{X_\reg}$, on (local) sections of ${\mathcal F}$ (with poles on reflection
hypertori). 

For a section $\beta$ of $(\cL^*)^w$ we define:
\begin{equation*}
\forall g\in W,\ (\rho(g)\beta)(x)=\beta(gx) 
\end{equation*}
(a section of $(\cL^*)^{gw}$); if $f$ is a section of ${\mathcal
O}_X$ then 
\begin{equation*}
\rho(f)\beta=f\beta;
\end{equation*}
and finally, for $v\in V$, 
\begin{equation*}
\rho(\partial_v)\beta=(\nabla_v+\sum_{g\in W}\langle 
F_{C,g}^{\cL^w},v\rangle)\beta.
\end{equation*}

\begin{proposition}\label{repre}
These formulas define a representation of $W\ltimes D_{X_\reg}$
on ${\mathcal F}|_{X_\reg}$. 
\end{proposition}

\begin{proof}
The only relations whose compatibility with $\rho$ needs to be
checked are $[\partial_v,\partial_u]=0$. 
This compatibility follows from statements (ii),(iv) of Lemma \ref{inv}. 
\end{proof} 

\begin{corollary}\label{repre1} The restriction of $\rho$ to $H_{c,X,W}\subset
W\ltimes D_{X_\reg}$ is a representation of $H_{c,X,W}$
on ${\mathcal F}$. 
\end{corollary}

\begin{proof}
We need to show that for any section $D$ of $H_{c,X,W}$,
$\rho(D)$ preserves holomorphic sections of ${\mathcal F}$.
Clearly, it is sufficient to check this for $D=D_{v,\phi}$, 
a Dunkl-Opdam operator. Obviously, we have 
$$
\rho(D_{v,\phi})|_{\cL^w}=\nabla_v+\sum_{(H,j)\in {\mathcal S}}C(H,j)
(\langle f_{H,j}^{\cL^w},v\rangle-\langle \phi_H,v\rangle
g_H^j).
$$
It is easy to see that each operator in parentheses preserves 
holomorphic sections, so the result follows. 
\end{proof} 

Note that the representation $\rho$ of $H_{c,X,W}$ belongs to
category ${\mathcal O}$, 
which is the category of representations of $H_{c,X,W}$ on
coherent sheaves on $X$. 

\section{Monodromy representation of orbifold Hecke algebras}

\subsection{Orbifold fundamental group and Hecke algebra}

The quotient $X/W$ is a complex orbifold. 
Thus for any $x\in X$ with trivial stabilizer, 
we can define the orbifold fundamental group $\pi_{1}^{\rm
orb}(X/W,x)$. It is the group consisting of the 
homotopy classes of paths on $X$ connecting $x$ and $gx$ for
$g\in W$, with multiplication defined by the rule: $\gamma_1\circ
\gamma_2$ is $\gamma_2$ followed by $g\gamma_1$, where $g$ is
such that $gx$ is the endpoint of $\gamma_2$. 
It is clear that the orbifold fundamental group of $X/W$ 
is naturally isomorphic to the semidirect product $W\ltimes \Gamma$. 

The braid group of $X/W$ is the orbifold fundamental group
$\pi_1^{\rm orb}(X_\reg/W,x)$. It can also be defined 
as $\pi_1(X'/W,x)$, where $X'$ is the set of all 
points of $X$ with trivial stabilizer. 

Now let $H$ be a reflection hypertorus.
Let $C_H$ be the conjugacy class in 
the braid group $\pi_1^{\rm orb}(X_\reg/W,x)$ corresponding
to a small circle going counterclockwise around 
the image of $H$ in $X/W$.  
Then we have the following result 
(see e.g. \cite{BMR}): 

\begin{proposition}\label{rel}
The group $\pi_1^{\rm orb}(X/W,x)=W\ltimes \Gamma$ is a quotient of the 
braid group $\pi_1^{\rm orb}(X_\reg/W,x)$ by 
the relations $T^{n_H}=1$ for all $T\in C_H$. 
\end{proposition}

Now for any conjugacy class of 
$H$, we introduce complex parameters
$\tau_{H,1},...,\tau_{H,n_H}$. 
The entire collection of these parameters 
will be denoted by $\tau$. 

The Hecke algebra of $(X,W)$, 
denoted ${\mathcal H}_\tau(X,W,x)$, is the quotient of 
the group algebra of the braid group, 
$\Bbb C[\pi_1^{\rm orb}(X_\reg/W,x)]$, by the relations
\begin{equation*}\label{poly}
\prod_{m=1}^{n_H} (T-e^{2\pi {\rm i}m/n_H}e^{\tau_{H,m}})=0,\ T\in C_H. 
\end{equation*}

This algebra is independent on the choice of $x$, so we will drop $x$ form the notation.

\begin{remark} 
It is known from \cite{E} that if $\tau$ is a formal, rather than
complex, parameter, 
then the algebra ${\mathcal H}_\tau(X,W)$ is a flat deformation
of the group algebra $\Bbb C[W\ltimes \Gamma]$. 
\end{remark}

\begin{example}\label{exa1}  If $W$ is a Weyl group, $V$ its reflection
representation, and $\Gamma=Q^\vee\oplus \eta Q^\vee$, where $Q^\vee$
is the dual root lattice of $W$, and $\eta\in \Bbb C^+$, then 
${\mathcal H}_\tau(X,W)$ is the double affine Hecke algebra 
(DAHA) of Cherednik, \cite{Ch2} (in the type $BC$ case one gets 
Sahi's 6-parameter version of the double affine Hecke algebra,
\cite{Sa}). 
\end{example}

\begin{example}\label{exa2} 
Let $W=S_n\ltimes (\Bbb Z/\ell Z)^n$, where $\ell=2,3,4,6$, 
$V=\Bbb C^n$, and $\Gamma=\Lambda^n$, where $\Lambda\subset \Bbb C$ is
a lattice invariant under $\Bbb Z/\ell Z$
(any lattice for $\ell=2$, triangular for $\ell=3,6$,
square for $\ell=4$). Then ${\mathcal H}_\tau(X,W)$ is the 
generalized double affine Hecke algebra of higher rank 
of type $D4,E6,E7,E8$, respectively, defined in \cite{EGO}.
We note that in the case $\ell=2$, this reproduces the 
BC case from Example \ref{exa1} (Sahi's algebra), and for $n=1$
these Hecke algebras were studied earlier in \cite{EOR}, in relation
to quantization of del Pezzo surfaces.  
\end{example} 

\subsection{The monodromy representation}

The representation $\rho$ defines a structure of a
$W$-equivariant holonomic $\mathcal O$-coherent D-module (i.e., a
$W$-equivariant local system) on the restriction of the vector bundle 
$\oplus_{w\in W}(\cL^*)^w$ to $X_\reg$. This local system yields
a monodromy representation $\pi_{\cL,\nabla}$ of the braid group
$\pi_1^{\rm orb}(X_\reg/W,x)$ (of dimension $|W|$). Since by Corollary \ref{repre1},
this local system is obtained by localization to $X_{\reg}$ 
of an ${\mathcal O}_X$-coherent 
$H_{c,X,W}$-module, by Proposition 3.4 of \cite{E}, 
the representation $\pi_\rho$ factors through the Hecke algebra 
${\mathcal H}_\tau(X,W)$, where $\tau$ is given by the formula
\begin{equation*}
\tau_{H,m}=-\frac{2\pi {\rm
i}}{n_H}\sum_{j=1}^{n_H-1}C(H,j)e^{-2\pi {\rm i}jm/n_H}.
\end{equation*}
Thus, for any collection of parameters $\tau_{H,j}$ with $\sum_j
\tau_{H,j}=0$ for all $H$, we have constructed a family of $|W|$-dimensional representations 
$\pi_{\cL,\nabla}$ of the Hecke algebra ${\mathcal H}_\tau(X,W)$,
parametrized by pairs $(\cL,\nabla)$; this family has $2\dim V$
parameters. 

Here is another version of the definition of the representation
$\pi_{\cL,\nabla}$ of ${\mathcal H}_\tau(X,W)$, which refers
directly to elliptic Dunkl operators and does not
mention elliptic Cherednik algebras. Consider the
system of differential-reflection equations
\begin{equation}\label{d=0}
\cD_{v,C}^{\cL,\nabla}\psi=0, v\in V. 
\end{equation}
Let ${\mathcal E}$ be the sheaf of solutions of this equation 
on $X'/W$ (sections of this sheaf over $\bar U=U/W$ are, by definition, solutions of
this system on $U$). Then ${\mathcal E}$ is a local system of
rank $|W|$, so it has a monodromy representation
$\xi_{\cL,\nabla}$. It is easy to see that
$\pi_{\cL,\nabla}=\xi_{\cL,\nabla}^*$. 

\begin{remark}
We could generalize the above construction by replacing
equations (\ref{d=0}) by the eigenvalue equations 
\begin{equation*}
\cD_{v,C}^{\cL,\nabla}\psi=\lambda(v)\psi, v\in V, 
\end{equation*}
where $\lambda\in V^*$, but this does not really give anything new
since it is equivalent to changing the connection $\nabla$. 
\end{remark}

\begin{example}
If $W$ is a Weyl group (Example \ref{exa1}), then the relation
$\sum_j \tau_{H,j}=0$ corresponds to the ``classical'' case 
of double affine Hecke algebras ($q=1$), in which case 
the DAHA is finitely generated over its center, and generically
over the spectrum of the center is an Azumaya algebra of rank
$|W|$, and the above construction yields  generic irreducible
representations of this algebra. In this case, such
representations can also be constructed by using classical analogs
of difference Dunkl-Cherednik operators \cite{Ch2}.    

If $W=S_n\ltimes (\Bbb Z/\ell \Bbb Z)^n$ (Example \ref{exa2}),
then we get (an open part of) the $2n$-parameter family of 
$|W|$-dimensional representations of generalized DAHA 
which was constructed in \cite{EGO} by another method. 
\end{example} 

In other cases of crystallographic reflection groups, however,  
the constructed family of representations appears to be new. 

\begin{conjecture} For any $W,V,\Gamma$, if $\sum_m
\tau_{H,m}=0$ for all $H$, then the Hecke algebra ${\mathcal H}_\tau(X,W)$ 
is finitely generated as a module over its center $Z_\tau(X,W)$, which is 
the algebra of functions on an irreducible affine algebraic
variety $M_\tau(X,W)$ of dimension $2\dim V$. Moreover, this
algebra is an Azumaya algebra of rank $|W|$ at the generic point
of $M_\tau(X,W)$, and the family of representations
$\pi_{\cL,\nabla}$ provides generic irreducible representations 
of ${\mathcal H}_\tau(X,W)$. 
\end{conjecture} 

This conjecture is known only in the case of Weyl groups (Example
\ref{exa1}), see \cite{Ch2}, and in the case $\dim V=1$, \cite{EOR}.
In particular, it is open in the case of Example \ref{exa2} 
for $\ell=3,4,6$.


\begin{thebibliography}{999}

\bibitem[BMR]{BMR} 
M. Brou\'e, G. Malle and R. Rouquier, Complex reflection groups,
braid groups, Hecke algebras, J. Reine Angew. Math. 500 (1998), 127-190.

\bibitem[BFV]{BFV} Buchstaber, V., Felder, G., Veselov, A., 
Elliptic Dunkl operators, root systems, and functional equations
Duke Math.J. 76 (1994) 885-911.

\bibitem[Ch1]{Ch1} Cherednik, I., Elliptic quantum many-body
problem and double affine Knizhnik-Zamolodchikov equation.
Comm. Math. Phys.  
169  (1995),  no. 2, 441--461.

\bibitem[Ch2]{Ch2} Cherednik, I. Double affine Hecke algebras,  
London Mathematical Society Lecture Note Series, 319, 
Cambridge University Press, Cambridge, 2005. 

\bibitem[Mu]{Mu}  D. Mumford, Abelian varieties, Oxford University Press, 1974.

\bibitem[DO]{DO}
    C. F. Dunkl, E. M. Opdam, Dunkl operators for complex
reflection groups, Proc. London Math. Soc. (3) 86 (2003), no. 1,
70-108.

\bibitem[E]{E} P. Etingof,
 Cherednik and Hecke algebras of varieties with a finite group action,
math.QA/0406499.

\bibitem[EGO]{EGO}
Etingof, P., Gan, W. L., Oblomkov, A., 
Generalized double affine Hecke algebras of higher rank.  
J. Reine Angew. Math.  600  (2006), 177--201.

\bibitem[EOR]{EOR} P. Etingof, A. Oblomkov, E. Rains,  
{\em Generalized double affine Hecke algebras of rank 1 and  
quantized del Pezzo surfaces}, arXiv:math/0406480,
to appear in Advances in Math.

\bibitem[GM]{GM} M. Geck, G. Malle, 
Reflection Groups, A Contribution to the Handbook of Algebra,
arXiv:math/0311012. 

\bibitem[La]{La} S. Lang, Abelian varieties, Springer-Verlag, New York, 1983.

\bibitem[Sa]{Sa} S. Sahi, Nonsymmetric Koornwinder polynomials
and duality, Ann. of Math. (2) 150 (1999), no. 1, 267--282.
\end{thebibliography}
\end{document}